\documentclass[11pt,bezier]{article}
\usepackage{amsmath}
\usepackage{amsfonts,amsthm,amssymb}
\usepackage{amsfonts}
\usepackage{graphics}
\usepackage{amssymb}
\newtheorem{lem}{Lemma}[section]
\newtheorem{thm}[lem]{Theorem}
\newtheorem{cor}[lem]{Corollary}

\theoremstyle{definition}

\begin{document}
\title{On the sizes of $k$-edge-maximal $r$-uniform hypergraphs
\footnote{The research is supported by NSFC (Nos. 11531011, 11771039, 11771443).}}
\author{Yingzhi Tian$^{a}$ \footnote{Corresponding author. E-mail: tianyzhxj@163.com (Y. Tian), xuliqiong@jmu.edu.cn (L. Xu), hjlai@math.wvu.edu (H. Lai), mjx@xju.edu.cn (J. Meng).}, Liqiong Xu$^{b}$, Hong-Jian Lai$^{c}$, Jixiang Meng$^{a}$ \\
{\small $^{a}$College of Mathematics and System Sciences, Xinjiang
University, Urumqi, Xinjiang 830046, PR China}\\
{\small $^{b}$School of Science, Jimei University, Xiamen, Fujian 361021, PR China}\\
{\small $^{c}$Department of Mathematics, West Virginia University,
Morgantown, WV 26506, USA}}

\date{}

\maketitle

\noindent{\bf Abstract } Let $H=(V,E)$ be a hypergraph, where $V$ is a  set of vertices and $E$ is a set of non-empty subsets of $V$ called edges. If all edges of $H$ have the same cardinality $r$, then $H$ is a $r$-uniform hypergraph; if $E$ consists of all $r$-subsets of $V$, then $H$ is a complete $r$-uniform hypergraph, denoted by $K_n^r$, where $n=|V|$. A hypergraph $H'=(V',E')$ is called a subhypergraph of $H=(V,E)$ if $V'\subseteq V$ and $E'\subseteq E$. A $r$-uniform hypergraph $H=(V,E)$ is $k$-edge-maximal if every subhypergraph of $H$ has edge-connectivity at most $k$, but for any edge $e\in E(K_n^r)\setminus E(H)$, $H+e$ contains at least one subhypergraph with edge-connectivity at least $k+1$.

Let $k$ and $r$ be integers with $k\geq2$ and $r\geq2$, and let $t=t(k,r)$ be the largest integer such that $(^{t-1}_{r-1})\leq k$. That is, $t$ is the integer satisfies $(^{t-1}_{r-1})\leq k<(^{t}_{r-1})$.
We prove that if $H$ is a $r$-uniform $k$-edge-maximal hypergraph such that $n=|V(H)|\geq t$, then

($i$) $|E(H)|\leq (^{t}_{r})+(n-t)k$, and this bound is best possible;

($ii$) $|E(H)|\geq (n-1)k -((t-1)k-(^{t}_{r}))\lfloor\frac{n}{t}\rfloor$, and this bound is best possible.

This extends former results in [8] and [6].


\section{Introduction}

In this paper, we consider finite simple graphs. For graph-theoretical terminologies and notation not defined here, we follow \cite{Bondy}.
For a graph $G$, we use $\kappa'(G)$ to denote the $edge$-$connectivity$ of $G$. The $complement$ of a graph $G$ is denoted by $G^c$. For $X\subseteq E(G^c)$, $G+X$ is the graph with vertex set $V(G)$ and edge set $E(G)\cup X$. We will use $G+e$ for $G+\{e\}$. The $floor$ of a real number $x$, denoted by $\lfloor x\rfloor$, is the greatest integer not larger than $x$; the $ceil$ of a real number $x$, denoted by $\lceil x\rceil$, is the least integer greater than or equal to $x$. For two integers $n$ and $k$, we define $(_k^n)=\frac{n!}{k!(n-k)!}$ when $k\leq n$ and $(_k^n)=0$ when $k>n$.

Given a graph $G$, Matula \cite{Matula} defined the $strength$ $\overline{\kappa}'(G)$ of $G$ as $max\{\kappa'(G'): G'\subseteq G\}$. For a positive integer $k$, the graph $G$ is $k$-$edge$-$maximal$ if $\overline{\kappa}'(G)\leq k$ but for any edge $e\in E(G^c)$, $\overline{\kappa}'(G+e)>k$. Mader \cite{Mader} and Lai \cite{Lai} proved the following results.

\begin{thm}
Let $k\geq1$ be an integer, and $G$ be a $k$-edge-maximal graph on $n>k+1$ vertices. Each of the following holds.

(i) (Mader \cite{Mader}) $|E(G)|\leq(n-k)k+(_2^k)$. Furthermore, this bound is best possible.

(ii) (Lai \cite{Lai}) $|E(G)|\geq (n-1)k-\lfloor\frac{n}{k+2}\rfloor(_2^k)$. Furthermore, this bound is best possible.
\end{thm}

In \cite{Anderson} and \cite{Lin}, $k$-edge-maximal digraphs are investigated, and the upper bound and the lower bound of the sizes of the $k$-edge-maximal digraphs are determined, respectively. Motivated by these results, we will study $k$-edge-maximal hypergraphs in this paper.

Let $H=(V,E)$ be a hypergraph, where $V$ is a finite set and $E$ is a set of non-empty subsets of $V$, called edges. An edge of cardinality 2 is just a graph edge. For a vertex $u\in V$ and an edge $e\in E$, we say $u$ is $incident$ $with$ $e$ or $e$ is $incident$ $with$ $u$ if $u\in e$.
If all edges of $H$ have the same cardinality $r$, then $H$ is a $r$-$uniform$ $hypergraph$; if $E$ consists of all $r$-subsets of $V$, then $H$ is a $complete$ $r$-$uniform$ $hypergraph$, denoted by $K_n^r$, where $n=|V|$.
For $n<r$, the complete $r$-uniform hypergraph $K_n^r$ is just the hypergraph with $n$ vertices and no edges.
The $complement$ of a $r$-uniform hypergraph $H=(V,E)$, denoted by $H^c$, is the $r$-uniform hypergraph with vertex set $V$ and edge set consisting of all $r$-subsets of $V$ not in $E$. A hypergraph $H'=(V',E')$ is called a $subhypergraph$ of $H=(V,E)$, denoted by $H'\subseteq H$, if $V'\subseteq V$ and $E'\subseteq E$.
For $X\subseteq E(H^c)$, $H+X$ is the hypergraph with vertex set $V(H)$ and edge set $E(H)\cup X$; for $X'\subseteq E(H)$, $H-X'$ is the hypergraph with vertex set $V(H)$ and edge set $E(H)\setminus X'$. We use $H+e$ for $H+\{e\}$ and $H-e'$ for $H-\{e'\}$ when $e\in E(H^c)$ and $e'\in E(H)$.
For $Y\subseteq V(H)$, we use $H[Y]$ to denote the hypergraph $induced$ by $Y$, where $V(H[Y])=Y$ and $E(H[Y])=\{e\in E(H): e\subseteq Y\}$. $H-Y$ is the hypergraph induced by $V(H)\setminus Y$.

For a hypergraph $H=(V,E)$ and two disjoint vertex subsets $X, Y\subseteq V$, let $E_H[X,Y]$ be the set of edges intersecting both $X$ and $Y$ and $d_H(X,Y)=|E_H[X,Y]|$. We use $E_H(X)$ and $d_H(X)$ for $E_H[X,V\setminus X]$ and $d_H(X,V\setminus X)$, respectively. If $X=\{u\}$, we use $E_H(u)$ and $d_H(u)$ for $E_H(\{u\})$ and $d_H(\{u\})$, respectively. We call $d_H(u)$ the $degree$ of $u$ in $H$.  The $minimum$ $degree$ $\delta(H)$ of $H$ is defined as $min\{d_H(u): u\in V\}$; the $maximum$ $degree$ $\Delta(H)$ of $H$ is defined as $max\{d_H(u): u\in V\}$. When $\delta(H)=\Delta(H)=k$, we call $H$ $k$-$regular$.

For a nonempty proper vertex subset $X$ of a hypergraph $H$, we call $E_H(X)$ an $edge$-$cut$ of $H$. The $edge$-$connectivity$ $\kappa'(H)$ of a hypergraph $H$ is $min\{d_H(X):\O\neq X\subsetneqq V(H)\}$. By definition, $\kappa'(H)\leq \delta(H)$. We call a hypergraph $H$ $k$-$edge$-$connected$ if $\kappa'(H)\geq k$. A hypergraph is connected if it is 1-edge-connected. A maximal connected subhypergraph of $H$ is called a $component$ of $H$. A $r$-uniform hypergraph $H=(V,E)$ is $k$-$edge$-$maximal$ if every subhypergraph of $H$ has edge-connectivity at most $k$, but for any edge $e\in E(H^c)$, $H+e$ contains at least one subhypergraph with edge-connectivity at least $k+1$. Since $\kappa'(K_n^r)=(_{r-1}^{n-1})$, we note that $H$ is a complete $r$-uniform hypergraph if $H$ is a $k$-edge-maximal $r$-uniform hypergraph such that $(_{r-1}^{n-1})\leq k$, where $n=|V(H)|$.
For results on the connectivity of hypergraphs, see [2,4,5] for references.

The main goal of this research is to determine, for given integers $n$, $k$ and $r$, the extremal sizes of a $k$-edge-maximal $r$-uniform hypergraph on $n$ vertices.
Section 2 below is devoted to the study of some properties of $k$-edge-maximal $r$-uniform hypergraphs. In section 3, we give the upper bound of the sizes of $k$-edge-maximal $r$-uniform hypergraphs and characterize these $k$-edge-maximal $r$-uniform hypergraphs attained this bound. We obtain the lower bound of the sizes of $k$-edge-maximal $r$-uniform hypergraphs and show that this bound is best possible in section 4.

\section{Properties of $k$-edge-maximal $r$-uniform hypergraphs}

For a $1$-edge-maximal $r$-uniform hypergraph $H$ with $n=|V(H)|$, we can verify that $\lceil\frac{n-1}{r-1}\rceil\leq|E(H)|\leq n-r+1$. If $H$ is the hypergraph with vertex set $V(H)=\{v_1,\cdots,v_n\}$ and edge set $E(H)=\{e_1,\cdots,e_{n-r+1}\}$, where $e_i=\{v_1,\cdots,v_{r-1},v_{r-1+i}\}$
for $i=1,\cdots,n-r+1$, then $H$ is  a $1$-edge-maximal $r$-uniform hypergraph $H$ with  $|E(H)|=n-r+1$. If $H$ is the hypergraph with vertex set $V(H)=\{v_1,\cdots,v_n\}$ and edge set $E(H)=\{e_1,\cdots,e_{s}\}$, where $s=\lceil\frac{n-1}{r-1}\rceil$, $e_i=\{v_{(i-1)(r-1)+1},\cdots,v_{i(r-1)},v_{n}\}$
for $i=1,\cdots,s-1$ and $e_s=\{v_{n-r+1}\cdots,v_{n-1},v_{n}\}$, then $H$ is  a $1$-edge-maximal $r$-uniform hypergraph $H$ with  $|E(H)|=\lceil\frac{n-1}{r-1}\rceil$.
Thus, from now on, we always assume $k\geq2$.

\noindent{\bf Definition 1.} For two integers $k$ and $r$  with $k,r\geq2$, define $t=t(k,r)$ to be the largest integer such that $(^{t-1}_{r-1})\leq k$. That is, $t$ is the integer satisfying $(^{t-1}_{r-1})\leq k<(^{t}_{r-1})$.

\begin{lem}
Let $H=(V,E)$ be a $k$-edge-maximal $r$-uniform hypergraph on $n$ vertices, where $k,r\geq2$.
Assume $n\geq t$ when $(^{t-1}_{r-1})= k$ and $n\geq t+1$ when $(^{t-1}_{r-1})<k$, where $t=t(k,r)$. Then $\kappa'(H)=\overline{\kappa}'(H)=k$.
\end{lem}

\noindent{\bf Proof.} Since $H$ is $k$-edge-maximal, we have $\kappa'(H)\leq\overline{\kappa}'(H)\leq k$. In order to complete the proof, we only need to show that $\kappa'(H)\geq k$.

Let $X$ be a minimum edge-cut of $H$, and let $H_1$ be a component of $H-X$ with minimum number of vertices and $H_2=H-V(H_1)$. Denote  $n_1=|V(H_1)|$ and $n_2=|V(H_2)|$. Thus we have $X=E_{H}[V(H_1), V(H_2)]$, $n=n_1+n_2$ and $n_1\leq n_2$. To prove the lemma, we consider the following two cases.

\noindent{\bf Case 1.} $E_{H^c}[V(H_1), V(H_2)]\neq\O$.

Pick an edge $e\in E_{H^c}[V(H_1), V(H_2)]$. Since $H$ is $k$-edge-maximal, we have $\overline{\kappa}'(H+e)>k$. Let $H'\subseteq H+e$
be a subhypergraph such that $\kappa'(H')\geq k+1$. By $\overline{\kappa}'(H)\leq k$, we have $e\in H'$. It follows that $(X\cup\{e\})\cap E(H')$ is an edge-cut of $H'$. Thus $|X|+1\geq |(X\cup\{e\})|\geq \kappa'(H')\geq k+1$, implying $|X|\geq k$. Thus $\kappa'(H)\geq k$.

\noindent{\bf Case 2.} $E_{H^c}[V(H_1), V(H_2)]=\O$.

Since $E_{H^c}[V(H_1), V(H_2)]=\O$, we know that $E_{H}[V(H_1), V(H_2)]$ consists of all $r$-subsets of $V(H)$ intersecting both $V(H_1)$ and $V(H_2)$. Thus
$$|E_{H}[V(H_1),V(H_2)]|=\sum_{s=1}^{r-1}(_s^{n_1})(_{r-s}^{n_2})=
(_r^{n})-(_r^{n_1})-(_r^{n_2}). $$
Let $g(x)=(_r^{x})+(_r^{n-x})$. It is routine to verify that $g(x)$ is a decreasing function when $1\leq x\leq n/2$. If $n_1\geq2$, then as $H$ is connected we have $r\leq n_1\leq n/2$. Thus
$$\kappa'(H)=|E_{H}[V(H_1),V(H_2)]|=
(_r^{n})-(_r^{n_1})-(_r^{n_2})\geq (_r^{n})-(_r^{2})-(_r^{n-2})>(_{r-1}^{n-1})\geq \delta(H), \eqno(1) $$
which contradicts to $\kappa'(H)\leq\delta(H)$. Thus, we assume $n_1=1$. Now we have
$$\kappa'(H)=|E_{H}[V(H_1),V(H_2)]|=
(_r^{n})-(_r^{n_1})-(_r^{n_2})= (_r^{n})-(_r^{1})-(_r^{n-1})=(_{r-1}^{n-1})\geq \delta(H), $$
which implies $\kappa'(H)=\delta(H)=(_{r-1}^{n-1})$ and so $H$ is a complete $r$-uniform hypergraph. Since $n\geq t$ when $(^{t-1}_{r-1})= k$ and $n\geq t+1$ when $(^{t-1}_{r-1})<k$, we have $\kappa'(H)=(_{r-1}^{n-1})\geq k$.
$\Box$

\begin{lem}
Suppose that $H=(V,E)$ is a $k$-edge-maximal $r$-uniform hypergraph, where $k,r\geq2$. Let $X\subseteq E(H)$ be a minimum edge-cut of $H$ and let $H_1$ be a union of some but not all components of $H-X$. Then $H_1$ is a $k$-edge-maximal $r$-uniform hypergraph.
\end{lem}

\noindent{\bf Proof.} If $H_1$ is complete, then $H_1$ is $k$-edge-maximal by definition. Thus assume $H_1$ is not complete.
For any edge $e\in E(H_1^c)\subseteq E(H^c)$, $H+e$ has a subhypergraph $H'$ with $\kappa'(H')\geq k+1$. Since $X$ is a minimum edge-cut of $H$, we have $|X|=\kappa'(H)\leq \overline{\kappa}'(H)\leq k$. Thus $X\cap E(H')=\O$. As $e\in E(H')\cap E(H_1^c)$, we conclude that $H'$ is a subhypergraph of $H_1+e$, and so $\overline{\kappa}'(H_1+e)\geq k+1$. Since $\overline{\kappa}'(H_1)\leq\overline{\kappa}'(H)\leq k$, it follows that $H_1$ is a $k$-edge-maximal $r$-uniform hypergraph.
$\Box$

\begin{lem}
Let $H=(V,E)$ be a $k$-edge-maximal $r$-uniform hypergraph on $n$ vertices, where $k,r\geq2$. Assume $n\geq t$ when $(^{t-1}_{r-1})= k$ and $n\geq t+1$ when $(^{t-1}_{r-1})<k$, where $t=t(k,r)$. Let $X\subseteq E(H)$ be a minimum edge-cut of $H$ and let $H_1$ be a union of some but not all components of $H-X$. If $r\leq|V(H_1)|\leq n-2$, then $|V(H_1)|\geq t$.  Moreover, if $H_1$ is complete, then $|V(H_1)|=t$; if $H_1$ is not complete, then $|V(H_1)|\geq t+1$.
\end{lem}

\noindent{\bf Proof.} By Lemmas 2.1 and 2.2, we have $|X|=\kappa'(H)=k$ and $H_1$ is a $k$-edge-maximal $r$-uniform hypergraph, respectively. If $H_1$ is not complete, then there is a subhypergraph $H_1'$ of $H_1+e$ such that $\kappa'(H_1')\geq k+1$ for any $e\in E(H_1^c)$. Since $(^{t-1}_{r-1})\leq k$ and $\delta(H_1')\geq \kappa'(H_1')\geq k+1$, we have
$|V(H_1)|\geq|V(H_1')|\geq t+1$.

Now we assume $H_1$ is a complete $r$-uniform hypergraph. Let $H_2=H-V(H_1)$. If $n_1=|V(H_1)|<t$, then, in order to ensure each vertex in $H_1$ has degree at least $k$ in $H$ (because $\delta(H)\geq \kappa'(H)=k$), we must have $n_1=t-1$ and $k=(^{t-1}_{r-1})$. Moreover, each vertex in $H_1$ is incident with exact $(^{t-2}_{r-2})$ edges in $E_H[H_1,H_2]$, and thus $d_H(u)=k$ for each $u\in V(H_1)$. By (1), there is an $e$ intersecting both $V(H_1)$ and $V(H_2)$ but $e\notin X$. Since $n_1\geq r$, there is a vertex $w\in V(H_1)$ such that $w$ is not incident with $e$. Then $d_{H+e}(w)=k$. This implies $w$ is not contained in a $(k+1)$-edge-connected subhypergraph of $H+e$. But then each vertex in $V(H_1)\setminus\{w\}$ has degree at most $k$ in $(H+e)-w$, and thus each vertex  in $V(H_1)\setminus\{w\}$ is not contained in a $(k+1)$-edge-connected subhypergraph of $H+e$. This illustrates that there is no $(k+1)$-edge-connected subhypergraph in $H+e$, a contradiction. Thus we have $n_1\geq t$. If $n_1> t$, then $\kappa'(H_1)=(^{n_1-1}_{r-1})\geq (^{t}_{r-1})>k$, contrary to $H$ is $k$-edge-maximal. Therefore, $n_1\leq t$,  and thus $n_1=t$ holds.
$\Box$

\section{The upper bound of the sizes of $k$-edge-maximal $r$-uniform hypergraphs}

\noindent{\bf Definition 2.} Let $n,k,r$ be integers such that $k,r\geq2$ and $n\geq t$, where $t=t(k,r)$. A hypergraph $H\in\mathcal{M}(n;k,r)$ if and only if it is constructed as follows:

($i$) Start from the complete hypergraph $H_0\cong K_t^r$;

($ii$) If $n-t=s=0$, then $H_s=H_0$. If $n-t=s\geq1$, then we construct, recursively, $H_i$ from $H_{i-1}$ by adding a new vertex $v_i$ and $k$ new edges containing $v_i$ and intersecting $V(H_{i-1})$ for $i=1,\cdots,s$;

($iii$) Set $H=H_s$.

It is known that $\kappa'(H)\leq\delta(H)$ holds for any hypergraph $H$. If $\kappa'(H)=\delta(H)$, then we say $H$ is $maximal$-$edge$-$connected$.
An edge-cut $X$ of $H$ is $peripheral$ if there exists a vertex $v$ such that $X=E_H(v)$. A hypergraph $H$ is $super$-$edge$-$connected$ if every minimum edge-cut of $H$ is peripheral. By definition, every super-edge-connected hypergraph is maximal-edge-connected.

\begin{lem}
Let $k$ and $r$ be integers with $k,r\geq2$. If $n\geq t$ when $(^{t-1}_{r-1})= k$ and $n\geq t+1$ when $(^{t-1}_{r-1})<k$, where $t=t(k,r)$, then for any $H\in \mathcal{M}(n;k,r)$, we have

(i) $\delta(H)=k$;

(ii) $H$ is super-edge-connected; and

(iii) $H$ is $k$-edge-maximal.
\end{lem}

\noindent{\bf Proof.} Let $H=H_s$, where $H_s$ is recursively constructed from $H_0,\cdots,H_{s-1}$ as in Definition 2. Then $V(H_s)=V(H_0)\cup\{v_1,\cdots,v_s\}$. We will prove this lemma by induction on $n$.

($i$) If $n=t$ and $(^{t-1}_{r-1})= k$, then $H\cong K_t^r$ and $\delta(H)=(^{t-1}_{r-1})=k$. If $n=t+1$ and $(^{t-1}_{r-1})<k$, then $H$ is obtained from $K_t^r$ by adding a new vertex $v_1$ and $k$ edges with cardinality $r$ such that each added edge is incident with $v_1$.  Let $k=(^{t-1}_{r-1})+i$. As $(^{t-1}_{r-1})<k<(^{t}_{r-1})$, we have $1\leq i\leq (^{t-1}_{r-2})-1$.
If there exists a vertex $u\in V(K_t^r)$ such that at most $i-1$ edges are incident with both $u$ and $v_1$ in $H$, then by $k=(^{t-1}_{r-1})+i$, we have $|E_H[\{v_1\}, V(H)\setminus\{u,v_1\}]|>(^{t-1}_{r-1})$. But this can not happen because $|V(H)\setminus\{u,v_1\}|=t-1$. Thus for any vertex $u\in V(K_t^r)$, there are at least $i$ edges incident with both $u$ and $v_1$ in $H$. This implies $d_H(v)\geq (^{t-1}_{r-1})+i=k$ for any $u\in V(K_t^r)$. As $d_H(v_1)=k$, we have $\delta(H)=k$.

Now we assume $n\geq t+1$ when $(^{t-1}_{r-1})= k$ and $n\geq t+2$ when $(^{t-1}_{r-1})<k$. Since $H=H_s$ is obtained from $H_{s-1}$ by adding a new vertex $v_s$ and $k$ edges with cardinality $r$ such that each added edge is incident with $v_s$, then by the induction assumption that $\delta(H_{s-1})=k$, we obtain $\delta(H)=\delta(H_s)=k$.

($ii$) If $n=t$ and $(^{t-1}_{r-1})= k$, then $H\cong K_t^r$ and $|E_H[X,V(H)\setminus X]|>\delta(H)=k$ for any $X\subseteq V(H)$ with $2\leq|X|\leq n-2$ by (1). Thus $H$ is super-edge-connected.

If $n=t+1$ and $(^{t-1}_{r-1})<k$, then $H$ is obtained from $K_t^r$ by adding a new vertex $v_1$ and $k$ edges with cardinality $r$ such that each added edge is incident with $v_1$. Let $k=(^{t-1}_{r-1})+i$. As $(^{t-1}_{r-1})<k<(^{t}_{r-1})$, we have $1\leq i\leq (^{t-1}_{r-2})-1$.
In order to prove that $H$ is super-edge-connected, we only need to verify that $d_H(X)>k$ for any  $X\subseteq V(H)\setminus\{v_1\}$ with $2\leq |X|\leq |V(H)|-2$. If $|X|\leq |V(H)|-3$, then $|E_{K_t^r}[X, V(K_t^r)\setminus X]|> (^{t-1}_{r-1})$ by (1). Since for any vertex $u\in V(K_t^r)$, there are at least $i$ edges incident with both $u$ and $v_1$ in $H$ (by the proof of $(i)$), we have $|E_{H}(X)\cap E_{H}(v_1)|\geq i$. Thus $d_H(X)=|E_{K_t^r}[X, V(K_t^r)\setminus X]|+|E_{H}(X)\cap E_{H}(v_1)|> (^{t-1}_{r-1})+i=k$. Assume $|X|= |V(H)|-2$ and $V(H)\setminus X=\{v_1,w\}$. If $r\geq3$, then $d_H(X)=|E_{K_t^r}[X, V(K_t^r)\setminus X]|+|E_{H}(X)\cap E_{H}(v_1)|=(^{t-1}_{r-1})+k>k$. If $r=2$, then $d_H(X)=|E_{K_t^r}[X, V(K_t^r)\setminus X]|+|E_{H}(X)\cap E_{H}(v_1)|\geq (^{t-1}_{r-1})+k-1>k$.

Now we assume $n\geq t+1$ when $(^{t-1}_{r-1})= k$ and $n\geq t+2$ when $(^{t-1}_{r-1})<k$. On the contrary, assume $H_s$ is not super-edge-connected. Then there is a minimum edge-cut $X=E_{H_s}[V(J_1), V(J_2)]$ of $H_s$ with $|X|\leq \delta(H_s)=k$, where $J_1$ is a component of $H_s-X$ and $J_2=H_s-V(J_1)$ with $min\{|V(J_1)|, |V(J_2)|\} \geq2$. Without loss of generality, assume $v_{s}\in V(J_1)$. If $E_{H_s}(v_s)\cap X\neq\O$, then as $X\neq E_{H_s}(v_s)$, $X-E_{H_s}(v_s)$ is an edge-cut of $H_{s-1}$, and so $\kappa'(H_{s-1})\leq |X-E_{H_s}(v_s)|< k$, contradicts to the induction assumption that $H_{s-1}$ is super-edge-connected. It follows that $E_{H_s}(v_s)\cap X=\O$ and so $X=E_{H_{s-1}}[V(J_1-v_s), V(J_2)]$ is an edge-cut of $H_{s-1}$. Since $H_{s-1}$ is super-edge-connected, we conclude that either $|V(J_1-v_s)|=1$ or $|V(J_2)|=1$. If $|V(J_2)|=1$, then it contradicts to $min\{|V(J_1)|, |V(J_2)|\} \geq2$. If $|V(J_1-v_s)|=1$, then $|V(J_1)|=2$, $r=2$ and $k=1$, contrary to $k\geq2$.

($iii$) If $n=t$ and $(^{t-1}_{r-1})= k$, then $H\cong K_t^r$ is $k$-edge-maximal by the definition.

If $n=t+1$ and $(^{t-1}_{r-1})<k$, let $k=(^{t-1}_{r-1})+i$. As $(^{t-1}_{r-1})<k<(^{t}_{r-1})$, we have $1\leq i\leq (^{t-1}_{r-2})-1$.
In order to prove that $H$ is $k$-edge-maximal, it suffices to verify that
$\overline{\kappa}'(H+e)\geq k+1$ for any $e\in E(H^c)$. By definition 2, $H+e$ is obtained from $K_t^r$ by adding a new vertex $v_1$ and $k+1$ edges with cardinality $r$ such that each added edge is incident with $v_1$.
If there exists a vertex $u\in V(K_t^r)$ such that at most $i$ edges are incident with both $u$ and $v_1$ in $H+e$, then by $k=(^{t-1}_{r-1})+i$, we have $|E_{H+e}[\{v_1\}, V(H)\setminus\{u,v_1\}]|>(^{t-1}_{r-1})$.  But this can not happen because $|V(H+e)\setminus\{u,v_1\}|=t-1$. Thus for any vertex $u\in V(K_t^r)$, there are at least $i+1$ edges incident with both $u$ and $v_1$ in $H+e$. This implies $d_{H+e}(u)\geq (^{t-1}_{r-1})+i+1=k+1$ for any $u\in V(K_t^r)$. By $d_{H+e}(v_1)=k+1$, we have $\delta(H+e)=k+1$.
For any edge-cut $W$ of $H+e$, if $W$ is peripheral, then $|W|\geq\delta(H+e)=k+1$. Suppose $W$ is not peripheral, and so $W-e$ is a non peripheral edge-cut of $H$. Since $H$ is super-edge-connected, $|W|\geq|W-e|\geq\delta(H)+1=k+1$. Thus $\overline{\kappa}'(H+e)\geq \kappa(H+e)\geq k+1$.

Now we assume $n\geq t+1$ when $(^{t-1}_{r-1})= k$ and $n\geq t+2$ when $(^{t-1}_{r-1})<k$. On the contrary, assume $H_s$ is not $k$-edge-maximal. Then there is an edge $e\in E(H_s^c)$
such that $\overline{\kappa}'(H_s+e)\leq k$. If $e\in E(H_{s-1}^c)$, then by induction assumption, $\overline{\kappa}'(H_{s-1}+e)\geq k+1$, a contradiction. Hence $e\notin E(H_{s-1}^c)$. Since $H_s$ is obtained from $H_{s-1}$ by adding a new vertex $v_s$ and $k$ edges incident with $v_s$, we have $e\in E_{H_{s}+e}(v_s)$.

Let $Y=E_{H_{s}+e}[V(F_1),V(F_2)]$ be a minimum edge-cut of $H_{s}+e$ with $|Y|\leq k$, where $F_1$ is a component of $(H_{s}+e)-X$ and $F_2=(H_{s}+e)-V(F_1)$. Since $H_s$ is super-edge-connected, we have $\kappa'(H_s)=\delta(H_s)=k$, and so $e\notin Y$ and $Y\neq E_{H_s}(v_s)$. This implies $Y\subseteq E(H_s)$. Without loss of generality, assume that $v_s\in V(F_1)$. By $H_{s-1}$ is super-edge-connected, we have $\kappa'(H_{s-1})=\delta(H_{s-1})=k$. If $Y\cap E_{H_s}(v_s)\neq\O$, then as
$Y\neq E_{H_s}(v_s)$, $Y-E_{H_s}(v_s)$ is an edge-cut of $H_{s-1}$. It follows that $\kappa'(H_{s-1})\leq |Y-E_{H_s}(v_s)|<k=\kappa'(H_{s-1})$, a contradiction. Hence we must have $Y\cap E_{H_s}(v_s)=\O$, and so $Y\subseteq E(H_s)-E_{H_s}(v_s)=E(H_{s-1})$. By $H_{s-1}$ is super-edge-connected, there exists a vertex $w\in V(H_{s-1})$ such that $Y=E_{H_{s-1}}(w)$. As $N_{H_s}(v_s)\cup\{v_s\}\subseteq V(F_1)$, we have $V(F_2)=\{w\}$.

Let $H'=H_s-w$. Then $e\in E((H')^c)$. If $w\in V(H_s)\setminus V(H_0)$, then $H'\in \mathcal{M}(n-1;k,r)$. If $w\in V(H_0)$, then by $d_{H_s}(w)=|Y|=k$, we have $d_{H_1}(w)=k$. By Definition 2,
there are exact $k-(^{t-1}_{r-1})$ edges containing  $\{w,v_1\}$ in $H_1$ and $|E_{H_1}[v_1,V(H_0)\setminus w]|=(^{t-1}_{r-1})$.
Thus the hypergraph induced by $(V(H_0)\setminus\{w\})\cup\{v_1\}$ in $H_s$ is complete, and so $H'\in \mathcal{M}(n-1;k,r)$. By induction assumption, $\overline{\kappa}'(H'+e)\geq k+1$, and so $\overline{\kappa}'(H_s+e)\geq\overline{\kappa}'(H'+e)\geq k+1$, contrary to  $\overline{\kappa}'(H_s+e)\leq k$.
$\Box$

\begin{thm}
Let $H$ be a $k$-edge-maximal $r$-uniform hypergraph on $n$ vertices, where $k,r\geq2$. If $n\geq t$, where $t=t(k,r)$, then each of the following holds.

(i) $|E(H)|\leq (^{t}_{r})+(n-t)k$.

(ii) $|E(H)|=(^{t}_{r})+(n-t)k$
if and only if $H\in \mathcal{M}(n;k,r)$.
\end{thm}

\noindent{\bf Proof.} By Definition 2, we have $|E(H)|=(^{t}_{r})+(n-t)k$ if $H\in \mathcal{M}(n;k,r)$.

We will prove the theorem by induction on $n$. If $n=t$, then by $H$ is $k$-edge-maximal and $(^{t-1}_{r-1})\leq k$, we have $H\cong K_t^r$. Thus  
$|E(H)|=(^{t}_{r})+(n-t)k$ and
$H\in \mathcal{M}(n;k,r)$.

Now suppose $n>t$.
We assume that if $t\leq n'<n$ and if $H'$ is a $k$-edge-maximal $r$-uniform hypergraph with $n'$ vertices, then $|E(H')|\leq (^{t}_{r})+(n'-t)k$ and $H'\in \mathcal{M}(n';k,r)$ if $|E(H')|=(^{t}_{r})+(n'-t)k$.

Let $X$ be a minimum edge-cut $H$. By Lemma 2.1, we have $|X|=k$. We consider two cases in the following.

\noindent{\bf Case 1.} There is a component, say $H_1$, of $H-X$ such that $|V(H_1)|=1$.

Let $H_2=H-V(H_1)$. By Lemma 2.2, $H_2$ is $k$-edge-maximal. Since $|V(H_2)|=n-1\geq t$, by induction assumption, we have $|E(H_2)|\leq (^{t}_{r})+(n-1-t)k$ and $H_2\in \mathcal{M}(n-1;k,r)$ if $|E(H_2)|=(^{t}_{r})+(n-1-t)k$. Thus $|E(H)|=|E(H_2)|+k\leq (^{t}_{r})+(n-t)k$. If $|E(H)|=(^{t}_{r})+(n-t)k$, then $|E(H_2)|=(^{t}_{r})+(n-1-t)k$ and $H_2\in \mathcal{M}(n-1;k,r)$. Thus, by $|V(H_1)|=1$ and $|X|=k$, we have 
$H\in \mathcal{M}(n;k,r)$ if $|E(H)|=(^{t}_{r})+(n-t)k$.

\noindent{\bf Case 2.} Each component of $H-X$ has at least two vertices.

Let $H_1$ be a component of $H-X$ and $H_2=H-V(H_1)$. By Lemma 2.2, both $H_1$
and $H_2$ are $k$-edge-maximal. Assume $n_1=|V(H_1)|$ and $n_2=|V(H_2)|$. Then $n_1+n_2=n$. Since each edge contains $r$ vertices, we have $n_1,n_2\geq r$. By Lemma 2.3, we have $n_1,n_2\geq t$. By induction assumption, we have $|E(H_i)|\leq (^{t}_{r})+(n_i-t)k$ and $H_i\in \mathcal{M}(n_i;k,r)$ if $|E(H_i)|=(^{t}_{r})+(n_i-t)k$ for $i\in \{1,2\}$. Thus

\ \ \ \ $|E(H)|=|E(H_1)|+|E(H_2)|+k$

\ \ \ \ \ \ \ \ \ \ \ \ \ \ $\leq(^{t}_{r})+(n_1-t)k+(^{t}_{r})+(n_2-t)k+k$

\ \ \ \ \ \ \ \ \ \ \ \ \ \ $=(^{t}_{r})+(n_1+n_2-t)k+(^{t}_{r})-(t-1)k$

\ \ \ \ \ \ \ \ \ \ \ \ \ \ $\leq(^{t}_{r})+(n_1+n_2-t)k+(^{t}_{r})-(t-1)(^{t-1}_{r-1})$

\ \ \ \ \ \ \ \ \ \ \ \ \ \ $=(^{t}_{r})+(n_1+n_2-t)k+(\frac{t}{r}-(t-1))(^{t-1}_{r-1})$

\ \ \ \ \ \ \ \ \ \ \ \ \ \ $\leq(^{t}_{r})+(n-t)k$.

If $|E(H)|=(^{t}_{r})+(n-t)k$, then $\frac{t}{r}-(t-1)=0$ and $k=(^{t-1}_{r-1})$, which imply $t=r=2$ and $k=1$, contrary to $k\geq2$.
Thus $|E(H)|<(^{t}_{r})+(n-t)k$ holds.
$\Box$

If $r=2$, then $H$ is a graph and $t=k+1$. Mader's \cite{Mader} result for the upper bound of the sizes of $k$-edge-maximal graphs is a corollary of Theorem 3.2.

\begin{cor} (Mader \cite{Mader})
Let $G$ be a $k$-edge-maximal graph with $n$ vertices, where $k\geq2$. If $n\geq k+1$, then we have  $|E(G)|\leq (^{k+1}_{2})+(n-k-1)k=(^{k}_{2})+(n-k)k$. Furthermore, $|E(G)|=(^{k}_{2})+(n-k)k$
if and only if $G\in \mathcal{M}(n;k,2)$.
\end{cor}

\section{The Lower bound of the sizes of $k$-edge-maximal $r$-uniform hypergraphs}

\begin{thm}
Let $H$ be a $k$-edge-maximal $r$-uniform hypergraph with $n$ vertices, where $k,r\geq2$. If $n\geq t$, where $t=t(k,r)$, then we have  $|E(H)|\geq (n-1)k -((t-1)k-(^{t}_{r}))\lfloor\frac{n}{t}\rfloor$.
\end{thm}

\noindent{\bf Proof.}
We will prove the theorem by induction on $n$. If $n=t$, then by $H$ is $k$-edge-maximal and $(^{t-1}_{r-1})\leq k$, we have $H\cong K_t^r$. Thus
$|E(H)|=(^{t}_{r})=(n-1)k -((t-1)k-(^{t}_{r}))\lfloor\frac{n}{t}\rfloor$.

Now suppose $n>t$.
We assume that if $t\leq n'<n$ and if $H'$ is a $k$-edge-maximal $r$-uniform hypergraph with $n'$ vertices, then $|E(H')|\geq (n'-1)k -((t-1)k-(^{t}_{r}))\lfloor\frac{n'}{t}\rfloor$.

Let $X$ be a minimum edge-cut $H$. By Lemma 2.1, we have $|X|=k$. We consider two cases in the following.

\noindent{\bf Case 1.} There is a component, say $H_1$, of $H-X$ such that $|V(H_1)|=1$.

Let $H_2=H-V(H_1)$. By Lemma 2.2, $H_2$ is $k$-edge-maximal. Since $|V(H_2)|=n-1\geq t$, by induction assumption, we have
$|E(H_2)|\geq (n-2)k -((t-1)k-(^{t}_{r}))\lfloor\frac{n-1}{t}\rfloor$. Thus 
 
\ \ \ \ $|E(H)|=|E(H_2)|+k$

\ \ \ \ \ \ \ \ \ \ \ \ \ \ $\geq(n-1)k -((t-1)k-(^{t}_{r}))\lfloor\frac{n-1}{t}\rfloor$

\ \ \ \ \ \ \ \ \ \ \ \ \ \ $\geq(n-1)k -((t-1)k-(^{t}_{r}))\lfloor\frac{n}{t}\rfloor$,

\noindent{the} last inequality holds because $(t-1)k-(^{t}_{r})\geq (t-1)(^{t-1}_{r-1})-\frac{t}{r}(^{t-1}_{r-1})\geq0$.

\noindent{\bf Case 2.} Each component of $H-X$ has at least two vertices.

Let $H_1$ be a component of $H-X$ and $H_2=H-V(H_1)$. By Lemma 2.2, both $H_1$
and $H_2$ are $k$-edge-maximal. Assume $n_1=|V(H_1)|$ and $n_2=|V(H_2)|$. Then $n_1+n_2=n$. Since each edge contains $r$ vertices, we have $n_1,n_2\geq r$. By Lemma 2.3, we have $n_1,n_2\geq t$. By induction assumption, we have $|E(H_i)|\geq(n_i-1)k -((t-1)k-(^{t}_{r}))\lfloor\frac{n_i}{t}\rfloor$ for $i\in \{1,2\}$. Thus

\ \ \ \ $|E(H)|=|E(H_1)|+|E(H_2)|+k$

\ \ \ \ \ \ \ \ \ \ \ \ \ \ $\geq(n_1-1)k -((t-1)k-(^{t}_{r}))\lfloor\frac{n_1}{t}\rfloor+(n_2-1)k -((t-1)k-(^{t}_{r}))\lfloor\frac{n_2}{t}\rfloor+k$

\ \ \ \ \ \ \ \ \ \ \ \ \ \ $=(n-1)k-((t-1)k-(^{t}_{r}))(\lfloor\frac{n_1}{t}\rfloor+
\lfloor\frac{n_2}{t}\rfloor)$

\ \ \ \ \ \ \ \ \ \ \ \ \ \ $\geq(n-1)k-((t-1)k-(^{t}_{r}))\lfloor\frac{n_1+n_2}{t}\rfloor$

\ \ \ \ \ \ \ \ \ \ \ \ \ \ $=(n-1)k-((t-1)k-(^{t}_{r}))\lfloor\frac{n}{t}\rfloor$.

The theorem is thus holds.
$\Box$

\noindent{\bf Definition 3.} Let $k, t, r$ be integers such that $t>r>2$, $k=(^{t-1}_{r-1})$ and $kr\geq 2t$. Assume $n=st$, where $s\geq2$. For any tree $T$ with $V(T)=\{v_1,\cdots,v_s\}$, we define a family of $r$-uniform hypergraphs
$\mathcal{N}(T)$ as follows. Firstly, we replace each $v_i$ by a complete $r$-uniform hypergraph $K_t^r(i)$ with $t$ vertices. Then whenever there is an edge $v_iv_j\in E(T)$, we add a set $E_{ij}$ of $k$ edges with cardinality $r$ such that ($i$) $e\subseteq V(K_t^r(i))\cup V(K_t^r(j))$, $e\cap V(K_t^r(i))\neq\O$ and $e\cap V(K_t^r(j))\neq\O$ for any $e\in E_{ij}$, and
($ii$) each vertex in $V(K_t^r(i))\cup V(K_t^r(j))$ is incident with some edge in $E_{ij}$ (we can do this because $kr\geq 2t$).

\begin{thm}
If $H\in \mathcal{N}(T)$, then $H$ is a $k$-edge-maximal $r$-uniform hypergraph.
\end{thm}

\noindent{\bf Proof.} By definition, $\overline{\kappa}'(H)\leq k$.  We will prove the theorem by induction on $s$. If $s=2$, then $|V(H)|=2t$ and $\delta(H)\geq (^{t-1}_{r-1})+1=k+1$. Since $K_t^r(1)$ and $K_t^r(2)$ are super-edge-connected and $\delta(H)\geq (^{t-1}_{r-1})+1=k+1$, each edge-cut of $H$ except for $E_H[V(K_t^r(1),V(K_t^r(2)]$ has cardinality at least $k+1$. For any $e\in E(H^c)$, we have $e\in E_{H^c}[V(K_t^r(1),V(K_t^r(2)]$. Thus every edge-cut of $H+e$ has cardinality at least $k+1$, that is, $\kappa'(H+e)\geq k+1$. This shows $\overline{\kappa}'(H+e)\geq \kappa'(H+e)\geq k+1$, and thus $H$ is $k$-edge-maximal.

Now suppose $s\geq 3$. We assume that each hypergraph constructed in Example 1 with less than $st$ vertices is $k$-edge-maximal. In the following, we will show that each $H$ in $\mathcal{N}(T)$ with $st$ vertices is also $k$-edge-maximal.

By contradiction, assume that there is an edge $e\in E(H^c)$ such that $\overline{\kappa}'(H+e)\leq k$. Let $E_{H+e}[X,V(H)\setminus X]$ be an edge-cut in $H+e$ with cardinality at most $k$. Since $K_t^r(i)$ is super-edge-connected for $1\leq i\leq s$ and $\delta(H)\geq k+1$, edge-cuts in $H$ with cardinality at most $k$ are these $E_{ij}$, where $v_iv_j\in E(T)$.
Thus $E_{H+e}[X,V(H)\setminus X]=E_{ij}$ for some $1\leq i,j\leq s$ with $v_iv_j\in E(T)$.  Then $e\in E_{H_i+e}(H_i+e)$, where $H_i$ is a component of $H-E_{ij}$. Since $H_i\in \mathcal{N}(T_i)$, where $T_i$ is a components of  $T-v_iv_j$, by induction assumption, $H_i+e$ contains a subhypergraph $H'$ with $\kappa'(H')\geq k+1$. But $H'$ is also a subhypergraph of $H+e$, contrary to $\overline{\kappa}'(H+e)\leq k$.
$\Box$

For any $H\in \mathcal{N}(T)$, we have $|E(H)|=(n-1)k -((t-1)k-(^{t}_{r}))\lfloor\frac{n}{t}\rfloor$. By Theorem 4.2, $H$ is $k$-edge-maximal. Thus, the lower bound given in Theorem 4.1 is best possible.

\end{document}